\documentclass{article}
\usepackage{latexsym, fullpage}

\usepackage{mathrsfs}
\usepackage[english]{babel}
\usepackage[utf8]{inputenc}
\usepackage{graphicx, color}
\usepackage[colorinlistoftodos]{todonotes}
\usepackage{fullpage}
\usepackage{graphics, color}
\usepackage{psfrag}
\usepackage{amsmath}
\usepackage{enumerate}
\usepackage{amsthm}
\usepackage{hyperref}
\usepackage{calrsfs}
\usepackage{tikz}
\usepackage{rotating}
\usepackage{amssymb}
\usepackage{amsmath}

\usepackage{authblk}
\usetikzlibrary{positioning,chains,fit,shapes,calc}

\newtheorem*{dfn}{Definition}

\newtheorem{theorem}{Theorem}[section]

\newtheorem{lemma}[theorem]{Lemma}

\theoremstyle{remark}

\title{The maximum spectral radius of planner graphs without the joint of $K_{2}$ and a linear forest}

\author{{Weilun Xu\thanks{Email: 225410010@fzu.edu.cn}}, {An Chang\thanks{Email: anchang@fzu.edu.cn}}
\\
\small{Center for Discrete Mathematics and Theoretical Computer Science\\ Fuzhou University\\ Fuzhou
Fujian, China}
 }

\date{ }

\begin{document}
\maketitle
\begin{center}
{\bf ABSTRACT}
\end{center}

Given a graph $F$, let $SPEX_P(n,F)$ be the set of graphs with the maximum spectral radius among all $F$-free $n$-vertex planner graph. In 2017, Tait and Tobin proved that for sufficiently $n$, $K_2+P_{n-2}$ is the unique graph with the maximum spectral radius
over all $n$-vertex planner graphs. In this paper, focusing on $SPEX_P(n,K_2+H)$ in which $H$ is a linear forest, we prove that $SPEX_P(n,K_2+H)=\{2K_1+C_{n-2}\}$ when $H\in \{pK_2,P_3,I_q\}$ $(p\geq1, q\geq 3)$, where $K_n$, $P_n$, $I_n$ are complete graph, path
and empty graph of order $n$, respectively. When $H$ contains a $P_4$, we prove that $2K_1+C_{n-2}\notin SPEX_P(n,K_2+H)$ and also provide a structural characterization of graphs in $SPEX_P(n,K_2+H)$.


\section{Introduction}
~

The Tur\'{a}n-type problems are undoubtedly one of the most extensively studied topics in extremal graph theory. Such a problem typically ask to maximize the number of edges in a graph which does not contain fixed forbidden subgraphs. Let $F$ be a graph. A graph
$G$ is called $F$-free if there is no subgraph of $G$ isomorphic to $F$. The earliest result on this issue was the well-known Mantel theorem\cite{Mantel} appeared in 1907, which says that every triangle free graph with $n$ vertices has at most
$\lfloor\frac{n^2}{4}\rfloor$  edges. In 1941, Tur\'{a}n\cite{Tur} extended this result to $K_{r+1}$-free graphs and proved that the $K_{r+1}$-free graphs with the maximum number of edges is the complete $r$-partite graph whose parts are of equal or almost equal
sizes, which is so called Tur\'{a}n graph. Given a graph $F$, the Tur\'{a}n number of $F$ is the maximum number of edges of an $F$-free graph with $n$ vertices. Nowadays the problem of determining the Tur\'{a}n number of a graph $F$ is known as the classical
Tur\'{a}n problem, which has received widespread attention and inspired a substantial research field of graph theory.
Meanwhile, there emerged many excellent results on this issue in which one of the most celebrated among them is ESS theorem, which is given by Erd\H{o}s, Stone and Simonovits.
 A graph $F$ is \emph{$k$-colorable} means that its vertex set can be colored with $k$ colors, such that every pair of adjacent vertices receives different colors. The \emph{chromatic number} of $F$, denoted by $\chi(F)$, is the minimum integer $k$ such that $F$
 is $k$-colorable. In \cite{Erdos46} and \cite{Erdos66}, Erd\H{o}s, Stone and Simonovits showed that every $F$-free graph with $n$ vertices and chromatic $\chi(F)$ has at most $(1-\frac{1}{\chi(F)-1}+o(1))\frac{n^2}{2}$ edges. The ESS theorem provides us a
 satisfied approximate answer to the Tur\'{a}n number of a graph $F$ when $\chi(F)\geq 3$ is known. Meanwhile, when $\chi(F)=2$, i.e., $F$ is a bipartite graph, we meet the so-called 'degenerate' problems which now is a major open problem on this issue. Such a
 problem is to determine the order of magnitude of Tur\'{a}n numbers for bipartite graphs.  For more details on this problem we refer the reader to the survey \cite{Furedi}.

For a given graph $G$ with vertex set $\{1,2,...,n\}$, the \emph{adjacent matrix} $A_G$ of $G$ is defined as follows: $$A_G(i,j)=\begin{cases}
1 & ij\in E(G),\\
0 & ij \notin E(G).
\end{cases}$$
The spectral radius of $G$, denoted by $\lambda(G)$, is the maximum eigenvalue of $A_G$.
One of the most well-known problems in spectral graph theory is the Brualdi–Solheid problem\cite{Brualdi}, which can state as the following

{\it{Problem 1}}. (Brualdi–Solheid problem) Given a set $\mathcal{G}$ of graphs, find a tight upper bound
for the spectral radius in $\mathcal{G}$ and characterize the extremal graphs.

When the graphs under consideration are $F$-free graphs, we have the following spectral version on the classical Tur\'{a}n problem

{\it{Problem 2}}. (Spectral Tur\'{a}n-type problem) What is the maximum spectral radius of a $F$-free graph
$G$ on $n$ vertices?

For a given graph $F$, we use $spex(n,F)$ to denote the {\emph{spectral Tur\'{a}n number}} of $F$, that is, the maximum spectral radius of an $n$-vertex $F$-free graph. As one of the pioneers of the extremal spectral graph theory, Nokiforov\cite{Nikiforov Kr}
proved in 2007 that $\lambda(G) \leq \lambda(T_r(n))$ when $G$ is an $n$-vertex $K_{r+1}$-free graph, where $T_r(n)$ is Tur\'{a}n graph on $n$-vertices. This result strengthened Tur\'{a}n's theorem because the spectral radius of any graph is no less than the
average degree of that graph.
 Additionally, Nokiforov also extend the result of Erd\H{o}s, Stone and Simonovits to a spectral version in \cite{Nikiforov ESS}.
For more about the spectral Tur\'{a}n problem, we refer the reader to the survey \cite{Li}.

A graph is planner if it can be drawn on the plane in such a way that its edges intersect only at their endpoints.~For a given planner graph $F$, the {\emph{planner Tur\'{a}n number} of $F$, denoted by $ex_{P}(n,F)$, is the maximum number of edges of an $n$-vertex
$F$-free planner graph. This topic was initiated by Dowden \cite{Dowden} in 2016. So far, the planner Tur\'{a}n number of some special graphs such as cycles \cite{Dowden} \cite{Ghosh} \cite{Gyori}\cite{Shi}, paths \cite{Lan19}\cite{Lan}, $\Theta$-graphs
\cite{Lan19DM}  have been determined.

It is naturally to consider the Tur\'{a}n-type spectral extremal problem for planner graph, i.e., to finde the maximum spectral radius over all $n$-vertex $F$-free planner graphs. We use $spex_P(n,F)$ to denote this value to be determined.
Instead of calculating $spex_P(n,F)$, we are interested in what graph can reach the maximum spectral radius among all $n$-vertex $F$-free planner graphs. Let $SPEX_{P}(n,F)$ be the set of all \emph{extremal graph}, that is, the $n$-vertex $F$-free planner graph
with spectral radius equals to $spex_P(n,F)$.

We use $C_{k}$, $I_{k}$ to denote a cycles on $k$ vertices, an empty graph with $k$ vertices, respectively. Let $G_1$ and $G_2$ be two graphs with disjoint vertex set. We use $G_1\cup G_2$ to denote the \emph{disjoint
union} of $G_1$ and $G_2$, i.e., $V(G_1\cup G_2)=V(G_1)\cup V(G_2)$ and $E(G_1\cup G_2)=E(G_1)\cup E(G_2)$. In particular, $kH$ is the disjoint union of $k$ copies of $H$. Further more, the \emph{join} of $G_1$ and $G_2$ is denoted by $G_1+G_2$, that is, $V(G_1+ G_2)=V(G_1)\cup V(G_2)$ and $E(G_1+ G_2)=E(G_1)\cup E(G_2)\cup
\{xy:x\in V(G_1)~\text{and}~y\in V(G_2)\}$.
Boots and Royle \cite{Boots}, and independently Cao and Vince \cite{Cao} conjectured that among all planner graph of order $n$, $K_2+P_{n-2}$ reaches the maximum spectral radius (Figure \ref{fig1}). In 2017, Tait and Tobin in their excellent paper\cite{Tait}
proved this conjecture for sufficiently large $n$.
The method utilized by Tait and Tobin in their proof has since been recognized as \emph{the second characteristic equation method}. This approach offers profound insights into the exploration of spectral  Tur\'{a}n problems.
 \begin{figure}[h]
 \begin{center}
  \includegraphics[scale=0.5]{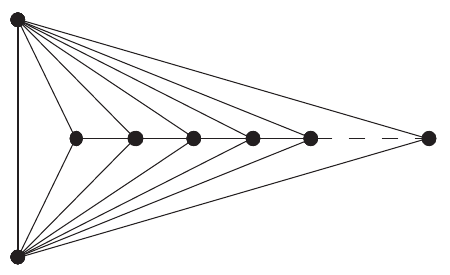}
\caption{{The graph $K_2+P_{n-2}$.\label{fig1}}}
\label{fig1}
 \end{center}
\end{figure}

An immediate application of the result obtained by Tait and Tobin \cite{Tait} is that, for $n$ large enough, $SPEX_P(n,F)=\{K_2+P_{n-2}\}$ when $F$ is not a subgraph of the graph $K_2+P_{n-2}$.  Consequently, when dealing with $SPEX_P(n, F)$ for sufficient large
$n$, it only becomes meaningful if $F$ is a subgraph of $K_2+P_{n-2}$. Recently, Fang, Lin and Shi \cite{Fang} proved that the unique extremal graph for the disjoint union of cycles is a join of a $K_2$ and a linear forest, where a \emph{linear forest} is a
disjoint union of paths and isolated vertices.

One can observe that a subgraph of the graph $K_2+P_{n-2}$ has chromatic number $4$ if and only if it contains a $K_4$. Furthermore, if such a subgraph is an induced subgraph, it must isomorphic to the join of a $K_2$ and a non-empty linear forest $H$. In the
light of these previous results,  this paper is concerned with the $SPEX_P(n,K_2+H)$ for $n$ large enough, where $H$ is a linear forest.

We first considered an arbitrary planner graph $F$ with chromatic number $4$ and provided a structural characterization of its spectral extremal graph.

\begin{theorem}\label{Lem:2case}
Let $F$ be a connected planner graph with chromatic number $4$. Suppose that $n$ is large enough and $G \in SPEX_P(n,F)$. Then one of the following holds:

(i)~$G\cong2K_1+C_{n-2}$.

(ii)~$K_2+I_{n-2}$ is a spanning subgraph of $G$.
\end{theorem}

Furthermore, we considered the two cases outlined in Theorem \ref{Lem:2case} in more detail.
We will demonstrate that $SPEX_P(n,K_2+H)$ can be determined by a parameter $\pi(H)$ of $H$, where $\pi(H)$ is defined as follows.
Let $H$ be a given linear forest. We use $ex^{F}(n,H)$ to denote the maximum number of edges of an $n$-vertex $H$-free linear forest and
$$\pi(H)=\lim_{n\rightarrow\infty}\frac{ex^{F}(n,H)}{n}.$$
For a given graph $H$, a graph $G$ is called \emph{$H$-maximal} if $G$ is $H$-free and adding any edge of the complement of $G$  will result in an $H$ in $G$. In particular, when $H$ is a linear forest, we say that a linear forest $H'$ on $n$ vertices is
$H$-maximal if $H'$ is $H$-free, and any linear forest on $n$ vertices that contains $H'$ as a proper subgraph is not $H$-free.
With the definition of $\pi(H)$, we are now in a position to present our main result as follows.

\begin{theorem}\label{Thm:main}
Suppose $H$ is a linear forest and $F=K_2+H$. For $n$ large enough,

(i)~If $\pi(H)< \frac{1}{2}$, then $SPEX_P(n,F)=\{2K_1+C_{n-2}\}$.

(ii)~If $\pi(H)> \frac{1}{2}$, then every $G\in SPEX_P(n,F)$ isomorphic to a $K_2+H'$, where $H'$ is an $H$-maximal linear forest.

\end{theorem}

According to Theorem \ref{Thm:main}, when $H$ is a linear forest and $\pi(H)\neq \frac{1}{2}$, then $\pi(H)$ can be regarded as a criterion to determine whether the the extremal graph of $F=K_2+H$ is $2K_1+C_{n-2}$ for $n$ large enough. The graph shown in
Figure \ref{fig2} is $2K_1+C_{8}$.
 \begin{figure}[h]
 \begin{center}
  \includegraphics[scale=0.5]{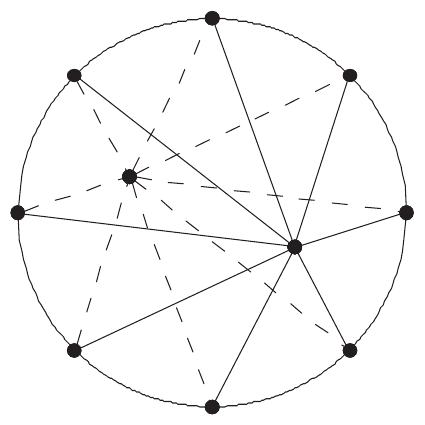}
\caption{{The graph $2K_1+C_{8}$.\label{fig2}}}
\label{fig1}
 \end{center}
\end{figure}

Starting from our main result, we can further investigate some interesting and meaningful cases. Let $pK_2$ be the disjoint union of $p$ edges. It is not difficult to demonstrate that $\pi(pK_2)<\frac{1}{2}$. Therefore, Theorem \ref{Thm:main} leads to the
following result.

\begin{theorem}\label{Cor:pK2}
For sufficiently large $n$ and $p\geq 1$, $SPEX_P(n,K_2+pK_2)=\{2K_1+C_{n-2}\}$.
\end{theorem}
We use $B_p$ to denote a \emph{book} on $p+2$ vertices, which is formally defined as $K_2+I_p$, where $I_p$ is an independent set on $p$ vertices. It is straightforward to observe that $B_p$ is a subgraph of $K_2+pK_2$. It is worth noting that $2K_1+C_{n-2}$ is
$B_p$-free for all $p\geq 3$. One immediate consequence of Theorem \ref{Cor:pK2} is the following result.
\begin{theorem}
For sufficiently large $n$ and $p\geq 3$, $SPEX_P(n,B_p)=\{2K_1+C_{n-2}\}$.
\end{theorem}

On the other hand, it is easy to verify that for $k\geq 4$, $\pi(P_k)>\frac{1}{2}$. Therefore, Theorem \ref{Thm:main} gives rise to the following two theorems.

\begin{theorem}\label{Thm:Pk}
For sufficiently large $n$, every $G\in SPEX_P(n,K_2+P_k)$ $(k\geq 4)$ isomorphic to a $K_2+H'$, where $H'$ is an  $H$-maximal linear forest.
\end{theorem}
\begin{theorem}
Let $H$ be a linear fores and $P_4\subseteq H$. Then for~$n$ large enough, every $G\in SPEX_P(n,K_2+H)$ isomorphic to a $K_2+H'$, where $H'$ is an $H$-maximal linear forest.
\end{theorem}
The remaining case of Theorem \ref{Thm:main} occurs when $\pi(H) = \frac{1}{2}$, indicating that the longest path present in graph $H$ is precisely a $P_3$. We begin by considering an instance of this, where $H=P_3$, and present the following theorem. In general, the graph $K_2+P_3$  is also referred to as  $K_5\setminus e$.
\begin{theorem}\label{thm:K5-}
For sufficiently large $n$, $SPEX_P(n,K_5\setminus e)=\{2K_1+C_{n-2}\}$.
\end{theorem}

The rest of this paper is organized as follows. In Section 2, we introduce the notations and terminologies which will be used frequently. A key lemma is also be given in this section. In Section 3, we mainly present the proof of our main results, following the order of Theorems \ref{Lem:2case}, \ref{Thm:main},
\ref{Cor:pK2}, \ref{Thm:Pk} and \ref{thm:K5-}.

\section{Preliminaries}
~

Let $G=(V,E)$ be a graph and $e(G)$ the number of edges in $G$. Suppose that $A$ and $B$ are two subsets of $V(G)$. We define $E(A,B)$ as the set of edges with one endpoint in $A$ and another in $B$. Furthermore, $e(A,B)$ denotes the cardinality of $E(A,B)$.
When $A$ and $B$ are the same, we simplify the notation to $E(A)$ and $e(A)$.  For any subset $A$ of $V(G)$, $G[A]$ represents the subgraph of $G$ induced by $A$. When referring to a specific vertex $v$, we use $N_G(v)$ to denote the neighbor of $v$ in $G$, and
$d_G(v)$ represents the number of neighbors of $v$, or $d_G(v)=|N_G(v)|$. Moreover, when discussing a subgraph $G_0$ of $G$, we define $N_{G_0}(v)$ as the intersection of $N_G(v)$ and $V(G_0)$, and $d_{G_0}(v)$ represents the number of neighbors of $v$ in $G_0$,
given by $d_{G_0}(v) = |N_{G_0}(v)|$.

For a given graph $G$, a \emph{minor} of $G$ is a graph that can be obtained from $G$ by means of a sequence of vertex and edge deletions and edge contractions. If $G$ does not contain $H$ as a minor, then $G$ is called $H$-minor free. The famous Kuratowski's
Theorem showed that a graph $G$ is planner if and only if $G$ is $\{K_5,K_{3,3}\}$-minor free. Another key property of planner graphs that we will use frequently is that a planner graph with $n$ vertices has at most $3n-6$ edges. If an $n$-vertex planner graph is
also bipartite, then the maximum number of edges is  $2n-4$.

In what follows, we always assume that the vertex set of an $n$-vertex graph $G$ considered in the paper is $\{1,2,...,n\}$. Let $\mathbf{v}=\{\mathbf{v}_1,\mathbf{v}_2,...,\mathbf{v}_n\}^{\intercal}$ be an $n$-dimensional vector over the $\mathbf{C}$. The
Euclidean norm ($l_2$-norm) and the max norm ($l_{\infty} $-norm) of $\mathbf{v}$ are
$\|\mathbf{v}\|=\sqrt{\mathbf{v}_1^2+\mathbf{v}_2^2+...+\mathbf{v}_n^2}$
 and $\|\mathbf{v}\|_{\infty}=\max{\{|\mathbf{v}_1|,|\mathbf{v}_2|,...,|\mathbf{v}_n|\}}$, respectively.
For a connected graph $G$ and positive constant $c$, by the well-known Perron-Frobenius Theorem in nonnegative matrix theory, there is a unique real eigenvector of $A_G$ in $\{\mathbf{v}\in\mathbb{R}^{|V(G)|}:\|\mathbf{v}\|_{\infty}=c\}$ corresponding to $\lambda(G)$, and it has strictly positive entries. We call such a vector \emph{Perron vector}. The
main approach we used to estimate $\lambda(G)$ is provided by Courant-Fischer Theorem, which states that $$\lambda(G)=\max_{\mathbf{v}\in \mathbb{R}^n}\frac{\mathbf{v}^\intercal A(G)\mathbf{v}}{\mathbf{v}^\intercal\mathbf{v}}.$$

Next, we define a commonly used graph transformation in the following text.
\begin{dfn}
Let $G$ be a graph and let $v \in V(G)$ be a vertex chosen from $G$. The graph $f_v^{V'}(G)$ is obtained by removing all edges incident with $v$ and connecting $v$ to every vertex in the subset $V' \subseteq V(G)$.

Moreover, for any edge $e \in E(G)$, the graph $f_{v,e}^{V'}(G)$ is defined as the graph derived from $G$ by removing all edges $\{uv \in E(G)\} \cup \{e\}$ and connecting $v$ to every vertex in $V' \subseteq V(G)$.
\end{dfn}

We need to show that the extremal graph contains a large bipartite subgraph $K_{2,t}$. Inspired by Tait and Tobin\cite{Tait}, we give the following lemma:
\begin{lemma}\label{Lem:fin}
Let $F$ be a connected planner graph with chromatic number $4$. For any $0<\epsilon\leq10^{-4}$ and $n$ large enough, let $G\in SPEX_{P}(n,F)$. Then \\
(i)~$G$ contains $K_{2,(1-100\epsilon)n}$ as a subgraph.\\
(ii)~If $\mathbf{v}$ is the Perron vector of $G$ with $\|\mathbf{v}\|_{\infty}=1$, then $\{\mathbf{v}_v:v\in V(G)\}\cap [900\epsilon,1-24\epsilon]=\emptyset$. Moreover, there exist exactly two vertices $x$ and $w$, which have maximum degree in the
$K_{2,(1-100\epsilon)n}$ such that $1=\mathbf{v}_x\geq \mathbf{v}_w>1-24\epsilon.$
\end{lemma}

{\it{Remark 2.}} In fact, Lemma \ref{Lem:fin} can be derived from the proof of Lemmas 8-12 in \cite{Tait}. It should be noted, however, that in the proof of Lemma \ref{Lem:fin}, when we perform a graph transformation, we must ensure that this graph transformation
does not alter the property of being $F$-free. This is the only significant difference between the proof of Lemma \ref{Lem:fin} and the proof of Lemma 8-12 in \cite{Tait}. Therefore, for the sake of completeness in our proof, we provide a detailed proof of Lemma
\ref{Lem:fin} and highlight the distinctions from the proof presented in \cite{Tait}.\\
{\bf\emph{Proof of Lemma \ref{Lem:fin}.}}
Let $F$ be a connected planner graph with chromatic number $4$.
Take an arbitrary $\epsilon \leq 10^{-4}$ and $G\in SPEX_{P}(n,F)$. Suppose $spex_P(n,F)=\lambda$ and $\mathbf{v}$ is the Perron vector of $G$ with $\|\mathbf{v}\|_{\infty}=1$. Let $x$ be a vertex such that $\mathbf{v}_x=1$.  Define
$$L=\{z\in V(G):\mathbf{v}_z>\epsilon\} ~\text{and}~ S=V(G) \setminus L.$$

By the definition of $L$, there is $$2(3n-6)\geq \sum_{z\in V(G)}d_G(z)\geq\sum_{z\in L}d_G(z)\geq \sum_{z\in V(G)}\lambda\mathbf{v}_z\geq |L|\epsilon\sqrt{2n-4},$$ and thus $|L|\leq \frac{3\sqrt{2n-4}}{\epsilon}\leq\epsilon n.$

Firstly, it is necessary to establish both the upper and lower bounds of $\lambda$. The following claim is derived from the facts that every $n$-vertex planer graph has at most $3n-6$ edges, and $K_2+I_{n-2}$ is an $F$-free planner graph. For a more comprehensive understanding of the proof of the subsequent claim, please refer to Lemma 8 of \cite{Tait}.\\
{\bf Claim 1.}\cite{Tait}
$$ \sqrt{2n-4}\leq\lambda\leq\sqrt{6n}.$$

~

Next, along the step of the proof of Lemma 9 in \cite{Tait}, we can estimate $\sum_{z\in L}\mathbf{v}_z$ and $\sum_{z\in S}\mathbf{v}_z$, respectively.\\
{\bf Claim 2.}\cite{Tait}
$$\sum_{z\in L}\mathbf{v}_z \leq\epsilon\sqrt{2n-4}+\frac{18}\epsilon$$ and $$\sum_{z\in S}\mathbf{v}_z\leq(1+3\epsilon)\sqrt{2n-4}.$$

~

Claim 2 demonstrates that in the summation $\sum_{y\sim x}\mathbf{v}_z=\lambda\geq \sqrt{2n-4}$ , the vertices belonging to the set $S$ account for the majority of the contribution. Furthermore, by the same calculation as in \cite{Tait} we can get that $$\sum_{\substack{y\in S\\ y\nsim x}}\mathbf{v}_y \leq 4\epsilon\sqrt{2n-4}+\frac{18}{\epsilon}.$$ Thus $|\{y\in S: y\nsim x\}|\mathbf{v}_{min}\leq 4\epsilon\sqrt{2n-4}+\frac{18}{\epsilon}$, where $\mathbf{v}_{min}=\min\{\mathbf{v}_y:y\in V(G)\}$. \\
{\bf Claim 3.}
$$\mathbf{v}_{min} \geq \frac{1}{\lambda} > \frac{1}{\sqrt{6n}}.$$

\proof
Suppose, to the contrary that there is a vertex $z$ with $\mathbf{v}_z < \frac{1}{\lambda}$. Then we have $z\nsim x$. For otherwise, $1=\mathbf{v}_x\leq \sum_{y\sim z}\mathbf{v}_y=\lambda\mathbf{v}_z<1$. A contradiction.

Let $G'=f_{z}^{\{x\}}(G)$. Then $G'$ is an $F$-free planner graph.
Otherwise, if $F\subseteq G'$, it must contain $z$. Since $n$ is sufficiently large, we have $d_{G'}(x)\geq d_G(x)\geq \lambda\geq\sqrt{2n-4}>|V(F)|$. There is a vertex $z'$ not in $V(F)$ and adjacent to $x$. Because $d_{G'}(z)=1$, $\{z'\}\cup(V(F)\cap
(V(G')\setminus z))$ contains a copy of $F$ in $G$. A contradiction.

Suppose that $\lambda(G')=\lambda'$. Now we have
\begin{align*}
\parallel \mathbf{v} \parallel^2\lambda' &\geq \mathbf{v}^\intercal A(G')\mathbf{v}\\
&=\mathbf{v}^\intercal A(G)\mathbf{v} -2\mathbf{v}_z(\sum_{y\sim z}\mathbf{v}_y)+2\mathbf{v}_z\mathbf{v}_x\\
&=\parallel \mathbf{v} \parallel^2\lambda +2\mathbf{v}_z(1-\lambda\mathbf{v}_z)\\
&>\parallel \mathbf{v} \parallel^2\lambda.
\end{align*}
A contradiction.

Now by $|\{y\in S: y\nsim x\}|\mathbf{v}_{min}\leq 4\epsilon\sqrt{2n-4}+\frac{18}{\epsilon}$ and Claim 3, we can give a lower bound of $d_{G}(x)$.\\
{\bf Claim 4.}\label{Lem:1}\cite{Tait}
There are at most $14\epsilon n$ vertices in $S$ non-adjacent to $x$. So $d_G(x)\geq (1-15)\epsilon n$.

~

Next, our objective is to locate a vertex $w\neq x$ with big degree. The proof of the subsequent claim follows a similar approach to the one utilized in the demonstration of Lemma 11 outlined in \cite{Tait}.\\
{\bf Claim 5.}\label{Lem:2}\cite{Tait}
There exists a vertex $w\neq x$ such that $\mathbf{v}_w > 1-24\epsilon$, and $$|\{y\in S:y\nsim w\}|\leq \sqrt{6n}(24\epsilon\sqrt{2n-4}-\frac{18}{\epsilon})\leq 84\epsilon n.$$

~

The remaining work is to bound $\mathbf{v}_z$ where $z\in V(G)\setminus\{x,w\}$. The proof is aligned with the approach outlined in Lemma 12 of \cite{Tait}.\\
{\bf Claim 6.}\label{Lem:3}\cite{Tait}
Let $v\in V(G)\setminus\{x,w\}$. Then $\mathbf{v}_v <900\epsilon<\frac{1}{10}$.

~

Clearly, Lemma \ref{Lem:fin} can be derived from Claim 4, 5 and 6.

\qed

\section{Proof of main results}
~
In this section, we always suppose that $n$ is a sufficiently large integer.

~

{\bf\noindent\emph{Proof of Theorem \ref{Lem:2case}.}}
Let $F$ be a connected planner graph with chromatic number $4$.
Suppose that $G\in SPEX_P(n,F)$. Let $\mathbf{v}$ be the Perron vector of $G$ with $\|\mathbf{v}\|_\infty=1$. Then by Lemma \ref{Lem:fin}, there exist two vertices $x$ and $w$ such that $1=\mathbf{v}_x\geq\mathbf{v}_w>1-24\epsilon$, where $\epsilon=10^{-4}$.
Denote $B=N(x)\cap N(w)$ and $A=V(G)\setminus (B\cup \{x,w\})$.

Note that $G[B]$ is $K_{1,3}$-free because $G$ contains no $K_{3,3}$. So $G[B]$ contains only cycles, paths and isolated vertices. Next, we are going to find a spanning subgraph of $G$ according to different cases.

 We will consider two cases based on whether there is a cycle in $G[B]$.\\
{\bf Case~1.}~There exists a cycle $C$ in $G[B]$.

We claim that $V(C)$ is a cut set such that $x$ and $w$ lies in different components of $G\setminus V(C)$. Otherwise, there exists a path $P$ with end vertices $x$ and $w$ such that $V(C)\cap V(P)=\emptyset$. Then $C$ is a $K_3$-minor and $P$ is a $K_2$-minor.
Since $\{x,w\}$ is complete to $C$, there is a $K_5$-minor. A contradiction. In particular, we have $V(C)=B$.

Let $D_x$ and $D_w$ be the  components of $G\setminus B$ that contain $x$ and $w$, respectively. Then we have $B\cup\{x,w\}\cup D_x \cup D_w =V(G)$. Otherwise, let $D$ be a component of $V(G)\setminus (B\cup\{x,w\}\cup D_x \cup D_w)$. We have $e(D)\leq 3|D|-6$
because $G[D]$ is a planner graph. Thus, there exists a vertex $z\in D$ has degree at most $5$ in $G[D]$. Moreover, $|N(z)\cap B|\leq 2$. (Otherwise, there is a $K_{3,3}$ in $G$.) Let $G'=f_{z}^{\{x\}}$. Clearly, $G'$ is an $F$-free planner graph. Suppose that
$\lambda(G')=\lambda'$. By Lemma \ref{Lem:fin}, we have
\begin{align*}
\parallel \mathbf{v} \parallel^2\lambda' &\geq \mathbf{v}^\intercal A(G')\mathbf{v}\\
&=\mathbf{v}^\intercal A(G)\mathbf{v} -2\mathbf{v}_z(\sum_{y\sim z}\mathbf{v}_y)+2\mathbf{v}_z\mathbf{v}_x\\
&=\mathbf{v}^\intercal A(G)\mathbf{v} -2\mathbf{v}_z(\sum_{\substack{y\in B\\y\sim z}}\mathbf{v}_y+\sum_{\substack{y\in D\\y\sim z}}\mathbf{v}_y)+2\mathbf{v}_z\mathbf{v}_x\\
&> \mathbf{v}^\intercal A(G)\mathbf{v} -2(\frac{2}{10}+\frac{5}{10})\mathbf{v}_z+2\mathbf{v}_z\mathbf{v}_x\\
&=\mathbf{v}^\intercal A(G)\mathbf{v} +\frac{3}{5}\mathbf{v}_z\\
&>\parallel \mathbf{v} \parallel^2\lambda.
\end{align*}
A contradiction.

Next, we are going to show that $D_x=\{x\}$ and $D_w=\{w\}$. Suppose to the contrary, there exists $z\in (D_x\cup D_w)\setminus\{x,w\}$. Let $p=|(D_x\cup D_w)\setminus\{x,w\}|\geq1$.\\
{\bf Claim.}~There exists a graph sequence $G=G_0,G_1,...,G_p=2K_1+C_{n-2}$ satisfies the following:

(i)~$G_p$ is an $F$-free planner graph.

(ii)~For any $0\leq i\leq p-1$, there exists an edge $u_iv_i\in E(C)\setminus(\cup_{j\leq i-1}\{u_jv_j\})$ and a vertex $z_i\in (D_x\cup D_w)\setminus(\cup_{j\leq i-1}\{z_j\}\cup\{x,w\} )$, such that $G_{i+1}=f^{\{x,w,u_i,v_i\}}_{z_i,u_iv_i}(G_{i})$.

(iii)~For any $1\leq i\leq p$, $\mathbf{v}^\intercal A(G_i)\mathbf{v}>\|\mathbf{v}\|^2\lambda$.

\proof

Note that the minimum subgraph with chromatic number $4$ of $G_p$  has at least $n-1$ vertices. We can assume that $n> 3|F|$. So $G_p$ is $F$-free. Thus (i) holds.

It is sufficient to show that for any $G_i$ ($0\leq i\leq p-1$), there exists $u_{i}v_{i}\in E(C)\setminus \cup_{j\leq i-1}\{u_jv_j\}$ and $z_{i}\in (D_x\cup D_w)\setminus (\cup_{j\leq i-1}\{z_j\}\cup\{x,w\})$ such that
$G_{i+1}=f^{\{x,w,u_i,v_i\}}_{z_i,u_iv_i}(G_{i})$ has $\mathbf{v}^\intercal A(G_i)\mathbf{v}>\|\mathbf{v}\|^2\lambda$.

For an edge $uv\in E(C)$, if $\mathbf{v}_u \leq \frac{60}{\lambda}$ and $\mathbf{v}_v \leq \frac{60}{\lambda}$, then $uv$ is called \emph{good}. Otherwise, $uv$ is called \emph{bad}. Let $L'=\{z\in V(G):\mathbf{v}_z>\frac{60}{\lambda}\}$. By the definition of $L'$ we have
$$2(3n-6)\geq \sum_{z\in V(G)}d_G(z)\geq\sum_{z\in L'}d_G(z)\geq \sum_{z\in L'}\lambda\mathbf{v}_z>60|L'|,$$
which implies $|L'|<\frac{1}{10}n.$ Note that $|E(C)|=|V(C)|=|N(x)\cap N(w)|\geq (1-100\epsilon)n$, thus $p\leq 100\epsilon n$.

Each vertex in $L'$ can contribute to at most two bad edges in $E(C)$. So there are at least $(1-100\epsilon)n-\frac{1}{5}n\geq p$ good edges. Take $p$ good edges in $E(C)$, and denote by $\{u_0v_0,u_1v_1,...,u_{p-1}v_{p-1}\}$.

 Note that for any $0\leq i\leq p-1$, $G[(D_x\cup D_w)\setminus (\cup_{j\leq i-1}\{z_j\})\cup\{x,w\}]$ is a planner graph. Consequently, there exists a vertex $z_{i}\in (D_x\cup D_w)\setminus (\cup_{j\leq i-1}\{z_j\})\cup\{x,w\}$ such that $z_i$ has at most $5$ neighbors in
 $(D_x\cup D_w)\setminus (\cup_{j\leq i-1}\{z_j\})\cup\{x,w\})$. Thus we get a vertex sequence $\{z_0,z_1,...,z_{p-1}\}$. In particular, each $z_i$ has at most two neighbors in $C$ and at most one neighbor in $\{x,w\}$.

 Recall that $\mathbf{v}_{x}\geq \mathbf{v}_{w}\geq 1-24\epsilon$ and $\frac{1}{\lambda}<\mathbf{v}_{z}<\frac{1}{10}$ for $z\notin \{x,w\}$. Now we have
\begin{align*}
 \mathbf{v}^\intercal A(G_{i+1})\mathbf{v}&\geq \mathbf{v}^\intercal A(G_i)\mathbf{v}+2\mathbf{v}_{z_i}(\mathbf{v}_{u_i}+\mathbf{v}_{v_i}+\mathbf{v}_{w}-\sum_{\substack{y\sim z_i\\y\neq x\\y\neq w}}\mathbf{v}_{y})-2\mathbf{v}_{u_i}\mathbf{v}_{v_i}\\
&>\mathbf{v}^\intercal A(G_i)\mathbf{v}+2\mathbf{v}_{z_i}(\mathbf{v}_{w}-\sum_{\substack{y\in C\\y\sim z_i}}\mathbf{v}_y-\sum_{\substack{y\in (D_x\cup D_w)\setminus (\cup_{j\leq i-1}\{z_j\}\cup\{x,w\})\\y\sim z_i}}\mathbf{v}_{y})-2\mathbf{v}_{u_i}\mathbf{v}_{v_i}\\
&>\mathbf{v}^\intercal A(G_i)\mathbf{v}+2\mathbf{v}_{z_i}(1-24\epsilon-\frac{2}{10}-\frac{5}{10})-\frac{7200}{\lambda^2}\\
&>\mathbf{v}^\intercal A(G_i)\mathbf{v}+\frac{2}{5\lambda}-\frac{7200}{\lambda^2}\\
&>\mathbf{v}^\intercal A(G_i)\mathbf{v}.
\end{align*}

Thus, $\|\mathbf{v}\|^2\lambda=\mathbf{v}^\intercal A(G_0)\mathbf{v}<\mathbf{v}^\intercal A(G_1)\mathbf{v}<\mathbf{v}^\intercal A(G_2)\mathbf{v}<...<\mathbf{v}^\intercal A(G_p)\mathbf{v}$. Therefore (ii) and (iii) hold, and we arrive the Claim.

From (i) and (iii) of above Claim, we can obtain a contradiction. This implies that $D_x=\{x\}$ and $D_w=\{w\}$. It follows that we have $G\cong 2K_1+C_{n-2}$.\\
{\bf Case 2.}~$G[B]$ is a linear forest.

We first show that $K_{2,n-2}\subseteq G$. Suppose to the contrary that $A\neq \emptyset$. Denote $|A|$ by $p$. Then there exists a vertex sequence $\{z_0,z_1,...,z_{p-1}\}$ such that $z_i$ has at most $5$ neighbors in $A\setminus (\cup_{j\leq i-1})\{z_j\}$. For
every $z\in A$, we have $|N(z)\cup B|\leq 2$. Now let $G_0=G$ and $G_{i+1}=f^{\{x,w\}}_{z_i}(G_{i})$ ($0\leq i\leq p-1$). Thus we have

\begin{align*}
 \mathbf{v}^\intercal A(G_{i+1})\mathbf{v}&\geq \mathbf{v}^\intercal A(G_i)\mathbf{v}+2\mathbf{v}_{z_i}(\mathbf{v}_{w}-\sum_{\substack{y\sim z_i\\y\neq x\\y\neq w}}\mathbf{v}_{y})\\
&>\mathbf{v}^\intercal A(G_i)\mathbf{v}+2\mathbf{v}_{z_i}(\mathbf{v}_{w}-\sum_{\substack{y\in C\\y\sim z_i}}\mathbf{v}_y-\sum_{\substack{A\setminus (\cup_{j\leq i-1})\{z_j\}}}\mathbf{v}_{y})\\
&>\mathbf{v}^\intercal A(G_i)\mathbf{v}+2\mathbf{v}_{z_i}(1-24\epsilon-\frac{2}{10}-\frac{5}{10})\\
&>\mathbf{v}^\intercal A(G_i)\mathbf{v}.
\end{align*}

Thus, $\|\mathbf{v}\|^2\lambda=\mathbf{v}^\intercal A(G_0)\mathbf{v}<\mathbf{v}^\intercal A(G_1)\mathbf{v}<\mathbf{v}^\intercal A(G_2)\mathbf{v}<...<\mathbf{v}^\intercal A(G_{p})\mathbf{v}$.

Note that $G_{p}=2K_1+(G[B]\cup I_{p})$ or $G_{p}=K_2+(G[B]\cup I_{p})$. So $G_{p}$ is a planner graph since $G[B]$ is a linear forest. Moreover, since $\lambda(G_p)>\lambda(G)$, we can infer that there is an $F$ in $G_p$. Then $G_p=K_2+(G[B]\cup I_{p})$ because $\chi(F)=4$. Let $F_1=V(F)\cap B$ and
$F_2=V(F)\setminus (F_1\cup\{x,w\})$. Note that $K_2+G[B]$ is both a subgraph of $G$ and $G_p$. Next, we are going to show that $K_2+G[B]$ contains a subgraph isomorphic to $F$. For a vertex $v\in F_2$, we have $d_F(v)\leq d_{G_p}(v) \leq 2$. Moreover, for a vertex $v\in B$, we
have $d_{G_p}(v) \geq 2$. Recall that $|B|\geq (1-100\epsilon)n > |F|$. So we can pick $F_3\subseteq B\setminus F_1$ with $|F_3|=|F_2|$. Then $F_1\cup F_3 \cup\{x,w\}$ contains a copy of $F$ in $G_p$, which means that $K_2+G[B]$  contains a copy of $F$. Thus, $G$ contains a
copy of $F$, which is a contradiction.

Note that if $x\nsim w$, then $G\subsetneq 2K_1+C_{n-2}$, contradiction with $G$ having maximum spectral radius. So $K_2+I_{n-2}$ is a spanning subgraph of $G$.\qed

~

{\bf\noindent\emph{Proof of Theorem \ref{Thm:main}.}} Let $\lambda(2K_1+C_{n-2})=\lambda$. Clearly, the Perron vector of the graph $2K_1+C_{n-2}$ takes the form of $(1,1,c,c,...,c)^\intercal$, where without loss of
generality, we suppose the components $1$ and $c$ correspond to the vertices of $2K_1+C_{n-2}$ with degrees $n-2$ and $4$, respectively. Additionally, using the system of eigen-equations
\begin{equation}
\left\{
\begin{aligned}\notag
(n-2)c &=\lambda &\\
2+2c &= c\lambda, &\\
\end{aligned}
\right.
\end{equation}

we obtain that $\lambda=\sqrt{2n-3}+1$ and $c=\frac{\lambda}{n-2}$.

(i)~Let $\pi(H)<\alpha< \frac{1}{2}$. Suppose that $G_1$ is the $F$-free graph with $n$ vertices and maximum spectral radius $\lambda_1$. If $G\neq 2K_1+C_{n-2}$, then by Theorem \ref{Lem:2case}, $G$ contains a $K_2+I_{n-2}$. Note that $K_{2,n-2}\subsetneq
G_{1}\subsetneq K_2+C_{n-2}$, which implies that $\sqrt{2n-4}<\lambda_1<\sqrt{2n-\frac{15}{4}}+\frac{3}{2}$, in which the last term is $\lambda(K_2+C_{n-2})$. In fact, to calculate $\lambda(K_2+C_{n-2})$, we consider the Perron vector $(1,1,c,c,\ldots,c)^{\intercal}$ of the graph $K_2+C_{n-2}$, where the components $c$ and $1$ of $\mathbf{v}$ correspond to the vertices of $K_2+C_{n-2}$ with degrees $4$ and $n-1$, respectively. By solving the following system of eigen-equations, we can calculate $\lambda(K_2+C_{n-2}) = \sqrt{2n - \frac{15}{4}} + \frac{3}{2}$
\begin{equation}
\left\{
\begin{aligned}\notag
1 + (n-2)c  &=\lambda(K_2+C_{n-2})&\\
2+2c &= c\lambda(K_2+C_{n-2}). &\\
\end{aligned}
\right.
\end{equation}

Let $x$ and $w$ be the vertices of $G_1$ with degree $n-2$. Suppose $\mathbf{v}$ is the Perron vector of $G_1$ with $\|\mathbf{v}\|_\infty=1$. Obviously, $\mathbf{v}_x=\mathbf{v}_w=1$. Recall that the common neighbor of $x$ and $w$, denoted by $B$, induces a linear forest.
So $G[B]$ has at most $\alpha n$ edges. Note that every $u\in B$ is adjacent to both $x$ and $w$. Thus $\mathbf{v}_u\geq \frac{2}{\lambda}$. Let $G_2=2K_1+C_{n-2}$ obtained from $G_1$ by removing $xw$ and extending $G_1[B]$ to a cycle. Then
\begin{align*}
 \mathbf{v}^\intercal A(G_{2})\mathbf{v}- \mathbf{v}^\intercal A(G_{1})\mathbf{v}&\geq 2((n-2-\alpha n)\frac{4}{\lambda^2} -1)\\
 &=2(\frac{4(1-\alpha)n+o(n)}{2n+o(n)}-1)>0.
\end{align*}
A contradiction. So (i) holds.

(ii)~Let $\pi(H)>\beta>\frac{1}{2}$. Then there exists an $H$-free linear forest $H'$ with $n-2$ vertices such that $|E(H')|>\beta n$. Let $G_1=K_2+H'$ and $G_2=2K_1+C_{n-2}$ obtained from $G_1$ by removing $xw$ and extend $H'$ to a $C_{n-2}$. Subsequently,
$G_2$ has spectral radius $\lambda=\sqrt{2n-3}+1$, and Perron vector $\mathbf{v}=(1,1,\frac{\lambda}{n-2},\frac{\lambda}{n-2},...,\frac{\lambda}{n-2})^\intercal$. Consequently,
\begin{align*}
 \mathbf{v}^\intercal A(G_{1})\mathbf{v}- \mathbf{v}^\intercal A(G_{2})\mathbf{v}&> 2(1+(\beta n-(n-2))\frac{\lambda^2}{(n-2)^2})\\
 &=2((\frac{2(\beta-1)n^2+o(n^2)}{n^2+o(n^2)}+1)>0.
\end{align*}
Thus $G\neq 2K_1+C_{n-2}$ . By Theorem \ref{Lem:2case}, $G$ has a spanning subgraph $K_2+I_{n-2}$. Let $H'=G[B]$. By the maximality of the spectral radius of $G$, $H'$ is an  $H$-maximal linear forest. So (ii) holds as desired.
\qed

~

{\bf\noindent{\emph{Proof of Theorem \ref{Cor:pK2}}}}. By Theorem \ref{Thm:main}, it is sufficient to show that $\pi(pK_2)=0<\frac{1}{2}$. Let $H'$ be a linear forest with $3p+1$ edges. Then we can embed a $pK_2$ into $H'$. We choose an arbitrary edge $e\in
E(H')$ to embed the first edge of $pK_2$. Suppose we have embedded $q$ edges. There are at most $2q$ edges in $E(H')$ adjacent to the embedded edges. So we still have $3(p-q)+1$ available edges. Consequently we can embed a $pK_2$ into $H'$. Thus we have
$$\pi(pK_2)=\lim_{n\rightarrow\infty}\frac{ex^{F}(n,pK_2)}{n}\leq\lim_{n\rightarrow\infty}\frac{3p+1}{n}=0.$$
\qed

~

{\bf\noindent{\emph{Proof of Theorem \ref{Thm:Pk}}}}.
By Theorem \ref{Thm:main}, it is sufficient to show that $\pi(K_2+P_k)> \frac{1}{2}$. Define $H_k(n)$ to be an $n$-vertex linear forest that is the disjoint union of $\lfloor\frac{n}{k-1}\rfloor$ copies of $P_{k-1}$ and a $P_{n-\lfloor\frac{n}{k-1}\rfloor(k-1)}$.
Clearly, $H_k(n)$ is $P_k$-free and
$$\pi(P_k)=\lim_{n\rightarrow\infty}\frac{ex^{F}(n,P_k)}{n}\geq \lim_{n\rightarrow\infty}\frac{|H_k(n)|}{n}=\frac{k-2}{k-1}>\frac{1}{2}.$$
\qed

~

{\bf\noindent{\emph{Proof of Theorem \ref{thm:K5-}}}}.
Theorem \ref{Lem:2case} shows that the extremal graph must be $2K_1+C_{n-2}$ or $K_2+H'$, where $H'$ is an $P_3$-maximal linear forest. Note that the unique $P_3$-maximal linear forest on $n$ vertices is the graph consists of some independent edges and at most one isolated point added. It suffices to demonstrate that
$\lambda(K_2+H')<\lambda(2K_1+C_{n-2})$. We first show that $\lambda(K_2+H')\leq \sqrt{2n-4}+1$. For notations simplicity, let $\lambda(K_2+H')=\mu$. \\
When $n$ is an odd integer, the Perron vector of the graph $K_2+H'$ has form $(1,1,c,c,...,c,d)^\intercal$, where the components $1$, $c$, and $d$ correspond to the vertices in $K_2+H'$ with degrees $n-1$, $3$ and $2$, respectively.
Furthermore, we can obtain the following system of eigen-equations:
\begin{equation}
\left\{
\begin{aligned}\notag
\mu d &= 2,&(1)\\
\mu c &= c+2,&(2)\\
\mu &= (n-3)c + d + 1.&(3)\\
\end{aligned}
\right.
\end{equation}
From the positivity of $c$, we infer that $c>d$ by equations (1) and (2). Substituting $c>d$ into the equation (3) leads to $(n-2)c+1>\mu$. Now, combining this inequality with the equation (2), we can deduce that $\mu <\sqrt{2n-4}+1$. \\
When $n$ is even, the system of eigen-equations simplifies to:
\begin{equation}
\left\{
\begin{aligned}\notag
\mu c &= c+2,&\\
\mu &= (n-2)c + 1.&\\
\end{aligned}
\right.
\end{equation}
A straightforward calculation leads to $\mu=\sqrt{2n-4}+1$. Recall that $\lambda(2K_1+C_{n-2})=\sqrt{2n-3}+1$ and thus it is easy to see that $\lambda(K_2+H')<\lambda(2K_1+C_{n-2})$.
\qed



\end{document}